\documentclass[12pt]{article}
\usepackage{graphicx}
\usepackage{amsmath,amsthm,amssymb,enumerate}
\usepackage{euscript,mathrsfs}
\usepackage[left=2cm,right=2cm,top=3.5cm,bottom=3.5cm]{geometry}
\usepackage{color}
\catcode`\@=11 \@addtoreset{equation}{section}

\catcode`\@=12

\allowdisplaybreaks

\newtheorem{Theorem}{Theorem}[section]
\newtheorem{Proposition}[Theorem]{Proposition}

\theoremstyle{definition}
\newtheorem{Definition}[Theorem]{Definition}

\newtheorem{Remark}[Theorem]{Remark}

\newcommand{\ep}{\varepsilon}
\newcommand{\vr}{\varrho}
\newcommand{\vt}{\vartheta}
\newcommand{\vc}[1]{{\bf #1}}
\newcommand{\vu}{\vc{u}}
\newcommand{\vm}{\vc{m}}
\newcommand{\Div}{{\rm div}_x}
\newcommand{\Grad}{\nabla_x}
\newcommand{\Ov}[1]{\overline{#1}}
\newcommand{\dx}{\,{\rm d}\vc{x}}
\newcommand{\dt}{\,{\rm d}t}
\newcommand{\dxdt}{\dx \,\dt}
\newcommand{\intO}[1]{\int_{\Omega}#1\dx\,}
\newcommand{\DC}{C^\infty_c}
\newcommand{\bfphi}{\boldsymbol{\varphi}}

\date{}

\begin{document}


\title{On oscillatory solutions to the complete Euler system}

\author{Eduard Feireisl
\thanks{The research of E.F.~leading to these results has received funding from the
European Research Council under the European Union's Seventh
Framework Programme (FP7/2007-2013)/ ERC Grant Agreement
320078. The Institute of Mathematics of the Academy of Sciences of
the Czech Republic is supported by RVO:67985840.}
\and Christian Klingenberg \thanks{C.K. and S.M. acknowledge the partial support by the ERC Grant Agreement 320078.}
\and Ond\v rej Kreml \thanks{The research of O.K. was supported by Neuron Fund for Support of Science under contract 18/2016 and by RVO:67985840.}
\and Simon Markfelder \footnotemark[2]}


\maketitle

\bigskip

\centerline{Institute of Mathematics of the Academy of Sciences of the Czech Republic}

\centerline{\v Zitn\' a 25, CZ-115 67 Praha 1, Czech Republic}
\medskip

\centerline{Department of Mathematics, W\"urzburg University}

\centerline{Emil-Fischer-Str. 40, 97074 W\"urzburg, Germany}

\bigskip

\begin{abstract}

The Euler system in fluid dynamics is a model
of a compressible inviscid fluid
incorporating the three basic physical principles: Conservation of mass, momentum, and energy.
We show that the Cauchy problem is basically ill--posed
for the $L^\infty$-initial data in the class of weak entropy solutions. As a consequence, there are infinitely many measure--valued solutions
for a vast set of initial data. Finally, using the concept of relative energy, we discuss a singular limit problem for the measure--valued solutions, where the Mach and Froude number are
proportional to a small parameter.

\end{abstract}

\tableofcontents

\section{Introduction}
\label{i}

The \emph{Euler system} is a model describing the time evolution of an inviscid fluid, the state of which is characterized by the mass density $\vr$, the (absolute) temperature $\vt$, and the macroscopic velocity field $\vu$. The problem can be written in the Eulerian coordinate system in the form:
\begin{align} \label{Ni1}
	\partial_t \vr + \Div (\vr \vu) &= 0, \\
	\label{Ni2}
	\partial_t (\vr \vu) + \Div (\vr \vu \otimes \vu) + \Grad p(\vr,\vt) &= \vr \Grad \Phi, \\
	\label{Ni3}
	\partial_t \left( \frac{1}{2} \vr |\vu|^2 + \vr e(\vr,\vt) \right) +
	\Div \left[ \left( \frac{1}{2} \vr |\vu|^2 + \vr e(\vr,\vt) + p(\vr, \vt) \right) \vu \right]  &= \vr \Grad \Phi\cdot \vu,
\end{align}
where $p$ is the pressure and $e$ the (specific) internal energy interrelated by Gibbs' equation
\begin{equation} \label{Ni4}
	\vt D s = D e + p D \left( \frac{1}{\vr} \right).
\end{equation}
The quantity $s = s(\vr, \vt)$ in (\ref{Ni4}) is the (specific) entropy. Any smooth solution of (\ref{Ni1}--\ref{Ni4}) satisfies the entropy balance
\begin{equation} \label{Ni5}
	\partial_t (\vr s(\vr, \vt))+ \Div( \vr s(\vr,\vt) \vu) = 0.
\end{equation}
The potential $\Phi$ describes the external force. In the absence of external forces, i.e. if $\Phi=0$, we call the Euler system (\ref{Ni1}-\ref{Ni3}) \emph{homogeneous}. The  non--homogeneous case was studied numerically e.g. in \cite{BEKR17,ChaKli15}.

Sometimes it is more convenient to rewrite the system (\ref{Ni1}--\ref{Ni3}) in the conservative variables
\[
\vr,\ \vm \equiv \vr \vu,\ E \equiv E_{\rm kin} + E_{\rm int},\ E_{\rm kin} \equiv \frac{1}{2} \vr|\vu|^2 = \frac{1}{2} \frac{|\vm|^2}{\vr},\
E_{\rm int} \equiv \vr e(\vr, \vt);
\]
\begin{align} \label{Ni6}
	\partial_t \vr + \Div \vm &= 0, \\
	\label{Ni7}
	\partial_t \vm + \Div \left( \frac{\vm \otimes \vm }{\vr} \right) + \Grad p(\vr,E) &= \vr \Grad \Phi, \\
	\label{Ni8}
	\partial_t E +
	\Div \left[ \Big( E + p(\vr, E) \Big) \frac{\vm}{\vr} \right] &= \Grad \Phi \cdot \vm.
\end{align}

In the real world applications, the fluid occupies a physical domain $\Omega \subset R^N$, $N=1,2,3$. For the sake of simplicity, we consider only the
natural \emph{impermeability} boundary conditions
\begin{equation} \label{Ni9}
	\vu \cdot \vc{n}|_{\partial \Omega} = 0.
\end{equation}
The original state of the fluid is determined by the initial conditions
\begin{equation} \label{Ni10}
	\vr(0, \cdot) = \vr_0,\ \vt(0, \cdot) = \vt_0,\ \vu(0, \cdot) =\vu_0.
\end{equation}

Solutions of (\ref{Ni1}--\ref{Ni3}), (\ref{Ni9}), (\ref{Ni10}) are known to develop singularities (shock waves) in a finite time for a rather generic class of the initial data,
see the classical monograph by Smoller \cite{Smoller67} or the more recent treatment by Benzoni--Gavage and Serre \cite{BenSer07}. To study the problem in the long run, the weak solutions must be considered. Unlike their classical counterparts, the weak solution may not satisfy the entropy balance (\ref{Ni5}) that must be relaxed to the \emph{inequality}
\begin{equation} \label{Ni11}
	\partial_t ( \vr s) + \Div (\vr s \vu) \geq 0.
\end{equation}
The weak solutions satisfying (\ref{Ni11}) are called \emph{entropy solutions} whereas relation (\ref{Ni11}) plays the role of a selection principle to identify
the physically relevant solution, cf. Dafermos  \cite{Dafermos73,Dafermos16}.

The entropy inequality (\ref{Ni11}) and its analogues for general hyperbolic systems proved to be efficient admissibility criterion in the simple 1-D setting.
In particular, Chen and Frid \cite{CheFri01,CheFri02} showed well--posedness of the 1-D Riemann problem for the system (\ref{Ni1}--\ref{Ni3}) in the class of
entropy solutions. The situation in the physically relevant multidimensional case $N=2,3$ turned out to be more delicate. In a series of papers,
De Lellis and Sz\' ekelyhidi  \cite{DelSze12,DelSze13,DelSze14} adapted the method of \emph{convex integration} to problems in fluid dynamics. Their effort
culminated by the complete proof of the celebrated Onsager's conjecture by Isett \cite{Isett17} and Buckmaster et al. \cite{BDSV17}. As a byproduct, new unfortunately mostly
negative results have been obtained concerning well--posedness of the \emph{isentropic} Euler system describing the motion of a compressible fluid with constant entropy,
see De Lellis and Sz\' ekelyhidi \cite{DelSze10}. Finally, Chiodaroli, De Lellis and Kreml \cite{ChiDelKre15} (see also
\cite{ChiKre14} for a more sophisticated example) showed that the isentropic Euler system is essentially ill--posed in the class of entropy (weak) solutions for Lipschitz initial data. Note that the
isentropic Euler system, where $p = p(\vr)$, contains only two unknowns, namely the density $\vr$ and the velocity $\vu$, whereas the entropy balance (\ref{Ni11}) is replaced
by the energy \emph{inequality}
\begin{align*}
	\partial_t \left( \frac{1}{2} \vr |\vu|^2 + P(\vr) \right) +
	\Div \bigg[ \bigg( \frac{1}{2} \vr |\vu|^2 + &P(\vr) + p(\vr) \bigg) \vu \bigg]\leq 0,\\
	&\mbox{where }P(\vr) = \vr \int_1^\vr \frac{p(z)}{z^2} \ {\rm d}z.
\end{align*}

This paper studies system (\ref{Ni1}--\ref{Ni3}) in the context of weak and more general measure--valued solutions. Combining a simple idea by
Luo, Xie, and Xin \cite{LuoXieXin16}  with a general framework developed in \cite{Feireisl16} we show that the Euler system (\ref{Ni1}--\ref{Ni3}), (\ref{Ni9}),
(\ref{Ni10}), supplemented
with the entropy inequality (\ref{Ni11}), admits infinitely many weak solutions for a large class of bounded initial data, see Section \ref{C}. This observation implies that the same problem admits genuine measure--valued solutions, meaning measure--valued solutions that do not coincide with a Dirac mass supported by a weak solution, 
see Section \ref{MVSo}. This kind of measure--valued solutions fits in the class \emph{dissipative measure--valued solutions} introduced in \cite{BreFei17}, but they are also the ``standard'' measure--valued solutions in the sense of Fjordholm, Mishra and Tadmor \cite{FjoMisTad16}.

Finally, in Section \ref{S}, using the relative energy inequality for the measure--valued solutions \cite{BreFei17}, we study the singular limit problem
for strongly stratified driven fluids described via the scaled system:
\begin{align} \label{Ni12} 
	\partial_t \vr + \Div (\vr \vu) &= 0, \\
	\label{Ni13}
	\partial_t (\vr \vu) + \Div (\vr \vu \otimes \vu) + \frac{1}{\ep^2} \Grad p(\vr,\vt) &= \frac{1}{\ep^2} \vr \Grad \Phi, \\
	\partial_t \left( \frac{1}{2} \vr |\vu|^2 + \frac{1}{\ep^2}\vr e(\vr,\vt) \right) 
	+\Div \left[ \left( \frac{1}{2} \vr |\vu|^2 + \frac{1}{\ep^2} \vr e(\vr,\vt) + \frac{1}{\ep^2} p(\vr, \vt)\right) \vu \right] &=\frac{1}{\ep^2} \vr \Grad \Phi \cdot \vu. \label{Ni14}
\end{align}
We consider the well--prepared initial data, where the density and the temperature are small perturbations of the \emph{isothermal} equilibrium state
\[
\vr_s = \vr_s(x), \ \Ov{\Theta} > 0
- \mbox{a positive constant,}\ \Grad p(\vr_s, \Ov{\Theta}) = \vr_s \Grad \Phi.
\]
We identify the limit system and show convergence of the dissipative measure--valued solutions of (\ref{Ni12}--\ref{Ni14}) for $\ep \to 0$.

\section{Well/ill--posedness}
\label{C}

We start by a definition of a weak solution to problem (\ref{Ni1}--\ref{Ni3}), (\ref{Ni9}), (\ref{Ni10}). To this end, it is more convenient
to use the formulation (\ref{Ni6}--\ref{Ni8}) based on conservative variables.

\begin{Definition} \label{defn:weaksol} {\bf [Entropy (weak) solution]} 
	
	Let $\Omega \subset R^N$, $N=1,2,3$ be a bounded domain. We say that a trio $[\vr, \vm, E]$ is an \emph{entropy (weak) solution} of problem
	(\ref{Ni6}--\ref{Ni10}) if
	\begin{itemize}
		\item $\vr \geq 0$ a.e. and the integral identity
		\begin{equation} \label{NC1}
			\int_0^T \intO{ \left[ \vr \partial_t \varphi + \vm \cdot \Grad \varphi \right] } \dt = - \intO{ \vr_0 \varphi (0, \cdot)}
		\end{equation}
		holds for any $\varphi \in \DC([0,T) \times R^N)$;
		\item $\vm = 0$ whenever $\vr = 0$,
		and the integral identity
		\begin{equation} \label{NC2}
			\begin{split}
				&\int_0^T \intO{ \left[ \vm \cdot \partial_t \bfphi + \frac{ \vm \otimes \vm }{\vr} : \Grad \bfphi + p(\vr, E) \Div \bfphi \right] } \dt \\
				&= - \intO{ \vm_0 \cdot \bfphi(0, \cdot) } - \int_0^T \intO{ \vr\Grad\Phi \cdot \bfphi \,} \dt 
			\end{split}
		\end{equation}
		holds for any $\bfphi \in \DC([0, T) \times \Ov{\Omega}; R^N)$, $\bfphi \cdot \vc{n}|_{\partial \Omega} = 0$;
		\item
		the integral identity
		\begin{equation} \label{NC3}
			\begin{split}
				&\int_0^T \intO{ \left[ E \partial_t \varphi + (E + p(\vr, E)) \frac{\vm}{\vr} \cdot \Grad \varphi \right] } \dt \\
				&= - \intO{ E_0 \varphi (0, \cdot) } - \int_0^T \intO{ \Grad\Phi\cdot \vm \,\varphi\, } \dt 
			\end{split}
		\end{equation}
		holds for any $\varphi \in \DC([0,T) \times R^N)$;
		\item
		the entropy inequality
		\begin{equation} \label{NC4}
			\begin{split}
				&\int_0^T \intO{ \left[ \vr Z \left(s(\vr, E) \right) \partial_t \varphi + Z\left(s(\vr, E) \right) \vm \cdot \Grad \varphi \right] } \dt \\
				&\leq - \intO{ \vr_0 Z(s(\vr_0, E_0)) \varphi (0, \cdot) }
			\end{split}
		\end{equation}
		holds for any $\varphi \in \DC([0,T) \times R^N)$, $\varphi \geq 0$, and any $Z \in C^2(R)$, $Z'(s) \geq 0$, $Z''(s) \leq 0$,
		$Z(s) \leq Z_\infty$ for any $s \in R$.
		
	\end{itemize}

\end{Definition}

\begin{Remark} \label{NR1}
	Here we have tacitly assumed that all integrals are well defined.
\end{Remark}

\begin{Remark}
	The entropy inequality is satisfied in the renormalized sense, similarly to Chen and Frid \cite{CheFri02}.
\end{Remark}

\subsection{Homogeneous Euler system}

In this section we consider the case $\Phi=0$. Furthermore, let $Q \subset R^N$, $N=2,3$ be a bounded domain, and $\vr > 0$, $p > 0$ positive constants. The nowadays standard result
based on convex integration (Chiodaroli \cite{Chiodaroli14}, Feireisl \cite[Theorem 13.6.1]{Feireisl16}) asserts that
there exists $\vm_0 \in L^\infty(Q; R^N)$ and a positive constant $\Lambda > 0$ such that the problem
\begin{align} \label{c4}
	\Div \vm &= 0,\\
	\label{c5}
	\partial_t \vm + \Div \left( \frac{ \vm \otimes \vm}{\vr} - \frac{1}{N} \frac{|\vm|^2}{\vr} \mathbb{I} \right) &= 0,\\
	\label{c6}
	\vm(0, \cdot) &= \vm_0
\end{align}
supplemented with the ``no--flux'' boundary conditions admits \emph{infinitely many} weak solutions $\vm$. In addition, these solutions satisfy
\begin{equation}
	\label{c7}
	E_{\rm kin} = \frac{1}{2} \frac{|\vm|^2}{\vr} = \Lambda - \frac{N}{2} p \quad \mbox{a.e. in}\ (0,T) \times Q,\
	E_{0,\rm kin} = \frac{1}{2} \frac{|\vm_0|^2}{\vr} = \Lambda - \frac{N}{2} p.
\end{equation}

By a weak solution with no--flux boundary conditions we mean a function
\begin{equation} \label{ClassM}
	\vm \in L^\infty((0,T) \times Q; R^N) \cap C_{\rm weak}([0,T]; L^2(Q; R^N))
\end{equation}
satisfying
\begin{align}
	\int_0^T \int_Q \vm \cdot \Grad \varphi \dxdt &= 0 \label{c8a} \\ \mbox{for any}\ \varphi \in &\DC ([0,T] \times R^N) \notag \\
	\int_0^T \int_Q \left[ \vm \cdot \partial_t \bfphi + \left( \frac{ \vm \otimes \vm}{\vr} - \frac{1}{N} \frac{|\vm|^2}{\vr} \mathbb{I} \right)
	: \Grad \bfphi \right] \dxdt &= - \int_Q \vm_0 \cdot \bfphi(0) \dx \label{c8b} \\
	\mbox{for any} \ \bfphi \in &\DC([0,T) \times R^N; R^N). \notag
\end{align}

\begin{Remark} \label{NR4}
	As a matter of fact, the method of convex integration provides infinitely many weak solutions for \emph{any} bounded initial momentum $\vm_0$ that
	may, and do in many cases, experience an initial jump of the kinetic energy. Here, the momentum $\vm_0$ gives rise to infinitely many weak solutions for which the initial energy is conserved.
\end{Remark}

\begin{Remark} \label{NR3}
	It is important that the momentum balance in (\ref{c8b}) holds for any smooth test function $\bfphi$, in particular, the normal trace
	of $\bfphi$ need not vanish on $\partial Q$. In particular, as observed by Luo, Xie, and Xin \cite{LuoXieXin16}, solutions may be defined piecewise
	on any finite union of mutually disjoint domains $Q_i$.
\end{Remark}

As $\vr > 0$ we can define a velocity field $\vu = \vm/\vr$. Moreover, we may assume $p = p(\vr, \vt)$ for another positive constant $\vt$.
Furthermore, using (\ref{c7}), (\ref{c8a}), (\ref{c8b}) we easily deduce
\begin{align} \label{c9a}
	\int_0^T \int_Q \left[ \vr \partial_t \varphi + \vr \vu \cdot \Grad \varphi \right] \dxdt &= -
	\int_Q \vr \varphi(0) \dx \\ \mbox{for any}\ \varphi \in &\DC ([0,T) \times R^N)\notag \\
	\int_0^T \int_Q \left[ \vr \vu \cdot \partial_t \bfphi + \left( { \vr \vu \otimes \vu} + p(\vr, \vt)\mathbb{I} - \frac{2}{N} \Lambda \mathbb{I} \right)
	: \Grad \bfphi \right] \dxdt &= - \int_Q \vm_0 \cdot \bfphi(0) \dx \label{c9b}\\
	\mbox{for any} \ \bfphi \in &\DC([0,T) \times R^N; R^N). \notag
\end{align}
Finally, we introduce the total energy,
\[
E = \frac{1}{2} \vr |\vu|^2 + \vr e(\vr,\vt) = E_{\rm kin} + \vr e(\vr, \vt),
\]
where, in accordance with (\ref{c7}), $E_{\rm kin}$ is a positive constant as long as $\Lambda > N/2 \ p$.
The quantities $E$ and $p$ being constant
in $(0,T) \times Q$, we easily deduce from (\ref{c8a}) the energy balance
\begin{equation} \label{c10}
	\begin{split}
		&\int_0^T \int_Q \left[
		\left( \frac{1}{2} \vr |\vu|^2 + \vr e(\vr, \vt) \right) \partial_t \varphi + \left[
		\left( \frac{1}{2} \vr |\vu|^2 + \vr e(\vr, \vt) + p(\vr, \vt) \right) \vu \cdot \Grad \varphi \right] \right] \dxdt \\ &= -
		\int_Q E \varphi(0) \dx \ \ \ \ \ \ \ \mbox{for any}\
		\varphi \in \DC([0,T) \times R^N),
	\end{split}
\end{equation}
together with the entropy balance
\begin{equation} \label{c11}
	\begin{split}
		\int_0^T \int_Q &\left[
		\vr Z(s(\vr, \vt))  \partial_t \varphi +  \vr Z(s(\vr, \vt))  \vu \cdot \Grad \varphi \right] \dxdt = -
		\int_Q \vr Z \left(s(\vr, \vt) \right) \varphi(0) \dx \\ &\mbox{for any}\
		Z\ \mbox{as in}\ \eqref{NC4} \mbox{ and any }
		\varphi \in \DC([0,T) \times R^N).
	\end{split}
\end{equation}

Finally, restricting the class of test functions $\bfphi$ in the momentum equation to those satisfying the boundary condition $\bfphi \cdot \vc{n}|_{\partial Q} = 0$
we may eliminate the $\Lambda$ dependent term in (\ref{c9b}). Thus we have shown that for any constant initial data $\vr = \vr_0 > 0$, $\vt = \vt_0 > 0$, there exists
$\vu_0 \in L^\infty$ such that the complete Euler system, supplemented with the impermeability condition $\vu \cdot \vc{n}|_{\partial Q} = 0$,
admits infinitely many weak entropy solutions.

As observed by Luo, Xie and Xin  \cite{LuoXieXin16}, the previous argument can be localized
and the result extended to piecewise--constant initial data $\vr_0$, $\vt_0$. It is also easy to observe that the parameter $\Lambda$ can be taken the same on
the whole domain $\Omega$. We say that a function $v \in L^\infty(\Omega)$ is piece--wise constant if
\[
\Omega = \cup_{i=1}^I \Ov{Q}_i, \ Q_i \subset R^N \ \mbox{domains},\ |\partial Q_i| = 0,\ \ Q_i \cap Q_j = \emptyset \ \mbox{for}\ i\ne j,\
v|_{Q_i} = v_i \ \mbox{- a constant.}
\]

\begin{Theorem} \label{Tc1}
	
	Let $\Omega \subset R^N$, $N=2,3$ be a bounded domain.
	Let the initial data $\vr_0 > 0$, $\vt_0 > 0$ be given piecewise constant functions in $L^\infty(\Omega)$.
	
	Then there exists $\vu_0 \in L^\infty(\Omega; R^N)$ such that the complete Euler system (\ref{NC1}--\ref{NC4}), with $\Phi=0$, admits
	infinitely many weak solutions originating from $[\vr_0, \vm_0, E_0]$,
	\[
	\vm_0 = \vr_0 \vu_0,\
	E_0 = \frac{1}{2} \vr_0 |\vu_0|^2 + \vr_0 e(\vr_0, \vt_0).
	\]
	In addition the renormalized entropy balance (\ref{NC4}) holds as equality, meaning the test functions need not be non-negative.

\end{Theorem}

We easily deduce from Theorem \ref{Tc1} that the set of initial densities and temperatures that gives rise to infinitely many solutions for certain
$\vu_0$ is dense in, say, $L^2(\Omega)$. Regularity properties of the corresponding momentum $\vm_0$ are not obvious. Note that the recent results
concerning Onsager's conjecture do not apply directly to problem (\ref{c4}--\ref{c6}) as the pressure term in (\ref{c5}) must be taken in a very particular form.

On the other hand, we report the following result proved by S. Markfelder and C. Klingenberg \cite{MarKli17}: There exists Lipschitz initial data for the isentropic Euler system which gives rise to infinitely many weak solutions that conserve energy. It is an easy observation that these energy conserving solutions are entropy solutions to the full Euler system, too (set $\vt$ such that $s=\text{const}.$). In other words there is Lipschitz data for the full Euler system (\ref{Ni1}--\ref{Ni3}) such that the Cauchy problem (\ref{Ni1}--\ref{Ni3}), (\ref{Ni10}) has infinitely many entropy solutions. These entropy solutions are similar to those in Definition \ref{defn:weaksol} with the only difference that the impermeability boundary conditions are not fulfilled.

\subsection{Driven fluids}

Next we consider driven fluids, i.e. $\Phi\neq 0$. Now we restrict ourselves to a 2-d bounded domain $\Omega \subset R^2$. Let again $\vr > 0$, $p > 0$ positive constants. As in the homogeneous case, the starting point is the convex integration result for the incompressible Euler equations described above: There exists $\vm_0\in L^\infty(\Omega;R^2)$ and a positive constant $\Lambda>0$ such that there are infinitely many solutions $\vm$ like in (\ref{ClassM}) fulfilling (\ref{c8a}), (\ref{c8b}) and
\begin{equation} \label{c7inhom}
	E_{\rm kin} = \frac{1}{2} \frac{|\vm|^2}{\vr} = \Lambda - p + \vr\,\Phi \quad \mbox{a.e. in}\ (0,T) \times \Omega,\
	E_{0,\rm kin} = \frac{1}{2} \frac{|\vm_0|^2}{\vr} = \Lambda - p + \vr\,\Phi.
\end{equation}
Note that this is a difference to (\ref{c7}) in the homogeneous case. Hence we deduce for the velocity field $\vu=\vm/\vr$ from (\ref{c8a}), (\ref{c8b}) and (\ref{c7inhom})
\begin{align} \label{c9ainhom}
	\int_0^T \intO{ \left[ \vr \partial_t \varphi + \vr \vu \cdot \Grad \varphi \right] }\dt &= -
	\intO{ \vr \varphi(0) } \\ \mbox{for any}\ \varphi \in &\DC ([0,T) \times R^2) \notag\\
	\int_0^T \intO{ \left[ \vr \vu \cdot \partial_t \bfphi + \left( { \vr \vu \otimes \vu} + p(\vr, \vt)\mathbb{I} - \Lambda \mathbb{I} - \vr\,\Phi \mathbb{I}\right)
		: \Grad \bfphi \right] }\dt &= - \intO{ \vm_0 \cdot \bfphi(0) }  \label{c9binhom}\\
	\mbox{for any} \ \bfphi \in &\DC([0,T) \times R^2; R^2), \notag
\end{align}
where again $\vt>0$ is such that $p=p(\vr,\vt)$.

Finally, we introduce the total energy,
\[
E = \frac{1}{2} \vr |\vu|^2 + \vr e(\vr,\vt) - \vr\,\Phi = E_{\rm kin} + \vr e(\vr, \vt) - \vr\,\Phi,
\]
where, in accordance with (\ref{c7inhom}), $E_{\rm kin}$ is positive as long as $\Lambda > p - \vr\,\Phi$.
The quantities $E$ and $p$ being constant
in $(0,T) \times \Omega$, we deduce from (\ref{c8a})
\begin{equation} \label{c10inhom}
	\begin{split}
		&\int_0^T \int_\Omega \bigg[
		\left( \frac{1}{2} \vr |\vu|^2 + \vr e(\vr, \vt) - \vr\,\Phi \right) \partial_t \varphi \\
		&\quad+ \left( \frac{1}{2} \vr |\vu|^2 + \vr e(\vr, \vt) + p(\vr, \vt) - \vr\,\Phi \right) \vu \cdot \Grad \varphi \bigg] \dx \dt \\ &= -
		\intO{ \left( \frac{1}{2} \vr |\vu|^2 + \vr e(\vr, \vt) - \vr\,\Phi \right)\varphi(0) } \ \ \ \ \ \ \ \mbox{for any}\
		\varphi \in \DC([0,T) \times R^2),
	\end{split}
\end{equation}
together with the entropy balance
\begin{equation} \label{c11inhom}
	\begin{split}
		\int_0^T &\intO{ \left[
			\vr Z(s(\vr, \vt))  \partial_t \varphi +  \vr Z(s(\vr, \vt))  \vu \cdot \Grad \varphi \right] } \dt = -
		\intO{ \vr Z \left(s(\vr, \vt) \right) \varphi(0) } \\ & \qquad \mbox{for any}\
		Z\ \mbox{as in}\ \eqref{NC4} \mbox{ and any }
		\varphi \in \DC([0,T) \times R^2).
	\end{split}
\end{equation}

Restricting to test functions $\bfphi$ in (\ref{c9binhom}) satisfying $\bfphi\cdot\vc{n}|_{\partial\Omega}$, we get rid of the term depending on $\Lambda$ and we can apply integration by parts in the term containing the potential $\Phi$. This procedure yields the correct momentum equation of Definition \ref{defn:weaksol}.

\begin{Remark}
	If we tried to consider piece--wise constant instead of constant densities as in the homogeneous case, we would get a problem. As above we can eliminate the term containing $\Lambda$ if we choose it the same on the whole domain. But the integration by parts of the term containing the potential $\Phi$ yields a boundary term that makes problems.
\end{Remark}

In addition it is easy to derive the energy equation of Definition \ref{defn:weaksol} from (\ref{c10inhom}) if we apply integration by parts to the term containing $\Phi$ and if we notice that $\vr\,\Phi$ does not depend on the time $t$.

Hence we proved the following:
\begin{Theorem} \label{Tc2}
	Let $\Omega \subset R^2$ be a bounded domain.
	Let the initial data $\vr_0 > 0$, $\vt_0 > 0$ be given constants.
	
	Then there exists $\vu_0 \in L^\infty(\Omega; R^2)$ such that the complete Euler system (\ref{NC1}--\ref{NC4}) admits
	infinitely many weak solutions originating from $[\vr_0, \vm_0, E_0]$,
	\[
	\vm_0 = \vr_0 \vu_0,\
	E_0 = \frac{1}{2} \vr_0 |\vu_0|^2 + \vr_0 e(\vr_0, \vt_0).
	\]
	In addition the renormalized entropy balance (\ref{NC4}) holds as equality, meaning the test functions need not be non-negative.
\end{Theorem}

\begin{Remark}
	Theorem \ref{Tc2} can be extended to 3-d cylindrical domains in a quite simple way as long as $\Phi=\Phi(y,z)$. Let $\Omega\subset\mathbb{R}^2$ be bounded and consider the domain $(0,1)\times\Omega$. For any constants $\rho_0,\vt_0>0$ - according to Theorem \ref{Tc2} - one can find $\vu_{0,h}\in L^\infty(\Omega,\mathbb{R}^2)$ such that the corresponding 2-d initial value problem has infinitely many entropy solutions $[\vr,\vm_h,E]$. Setting $\vm=(0,\vm_h)$ and keeping in mind that the first component of $\Grad\Phi$ is zero, it is easy to show that $[\vr,\vm,E]$ is an entropy solution to the 3-d initial value problem with initial data $\vr_0,\vt_0,\vu_0:=(0,\vu_{0,h})$.
\end{Remark}

\begin{Remark}
	The results of Section \ref{C}, i.e. Theorems \ref{Tc1} and \ref{Tc2}, extend to unbounded domains with obvious modifications.
\end{Remark}

\section{Measure--valued solutions}
\label{MVSo}

The so--called measure--valued solutions have been introduced in fluid dynamics in the pioneering paper by DiPerna and Majda \cite{DipMaj87}, and revisited
recently in the context of numerical analysis by Fjordholm, Mishra, and Tadmor \cite{FjoMisTad12,FjoMisTad16}. The measure--valued solutions capture possible
oscillations and concentrations that may appear in families of (approximate) solutions to the Euler system. Roughly speaking, the exact values of non--linear
compositions are replaced by their expectations with respect to a probability measure (Young measure). Theorem \ref{Tc1} above indicates that the measure--valued solutions may be indeed relevant in the context of inviscid fluids. Indeed any weak solution $[\vr, \vt, \vu]$ of problem (\ref{NC1}--\ref{NC4}) can be interpreted as a measure--valued solution represented by the Dirac measure
\[
\delta_{\vr(t,x), \vt(t,x), \vu(t,x)}
\]
supported by the weak solution. In view of Theorem \ref{Tc1}, the Euler system (\ref{NC1}--\ref{NC4}) admits infinitely many solutions
$[\vr_i, \vt_i, \vu_i]$, $i \in I$
for certain initial data. In particular, the family of parameterized measures
\[
\sum_{i=1}^{I_0} \lambda_i \delta_{\vr_i(t,x), \vt_i(t,x), \vu_i(t,x)}, \ \sum_{i=1}^{I_0} \lambda_i = 1
\]
represents a (genuine) measure--valued solution of the same problem (cf. Definition \ref{DMV1} below).

A proper definition of measure--valued solution to the full Euler system is a delicate issue. As uniform $L^\infty$ {\it a priori} bounds are not known and possibly not available, the effect of possible concentrations must be taken into account. Unfortunately, the $L^p-$framework developed in \cite{FjoMisTad16}
does not apply either as the fluxes are not dominated by the conserved quantities. To avoid this difficulty, we follow the approach \cite{BreFei17}, where
the energy equation (\ref{NC3}) is ``integrated'' with respect to the space variable and the entropy balance kept in its renormalized form (\ref{NC4}).
In such a way, we eliminate the difficulties connected with the lack of {\it a priori} bounds that would control the convective terms in (\ref{NC3}), (\ref{NC4}).

\subsection{Dissipative measure--valued solutions}

Following \cite{BreFei17} we define the measure--valued solutions to the Euler system in terms of variables
\[
\vr, \ \vm, \ \mbox{and}\ E_{\rm int} = \vr e(\vr, \vt)
\]
assuming that $\partial_\vt e(\vr, \vt) > 0$. Accordingly, we consider the phase space
\[
\mathcal{F} = \left\{ [\vr, \vm, E_{\rm int}] \ \Big| \ \vr \geq 0, \ \vm \in R^N, \ E_{\rm int} \geq 0 \right\}.
\]

\begin{Definition} \label{DMV1} {\bf [Dissipative measure--valued solution]} 
	
	A parameterized family of probability measures $\{ Y_{t,x} \}_{(t,x) \in (0,T) \times \Omega}$,
	\[
	(t,x) \mapsto Y_{t,x} \in L^\infty_{\rm weak-(*)} ((0,T) \times \Omega; \mathcal{P}(\mathcal{F})),\
	\]
	and a non-negative function $\mathcal{D} \in L^\infty(0,T)$ called dissipation defect
	represent a \emph{dissipative measure--valued solution} of the Euler system
	(\ref{Ni1}--\ref{Ni3}), (\ref{Ni9}), (\ref{Ni11})
	with the initial data $Y_{0,x}$ if:
	\begin{itemize}
		\item 
		\begin{equation} \label{MV22}
			\left[ \intO{
				\left< Y_{t,x}; \vr \right> \varphi } \right]_{t = 0}^{t = \tau} =
			\int_0^\tau \intO{ \left[ \left< Y_{t,x}; \vr \right> \partial_t \varphi + \left< Y_{t,x}; \vm \right> \cdot \Grad \varphi \right] } \dt
		\end{equation}
		for a.a. $\tau \in (0,T)$ and for any $\varphi \in C^1 ([0,T] \times R^N)$;
		
		\item
		\begin{equation} \label{MV23}
			\begin{split}
				&\left[ \intO{ \left< Y_{t,x}; \vm \right> \cdot \bfphi }
				\right]_{t = 0}^{t = \tau} \\ &= \int_0^\tau \int_\Omega \bigg[ \left< Y_{t,x}; \vm \right> \cdot \bfphi + \left< Y_{t,x}; \frac{\vm \otimes \vm }{\vr}  \right> : \Grad \bfphi + \left< Y_{t,x} ; p(\vr,E_{\rm int}) \right> \Div \bfphi \bigg] \dx \dt\\
				&\quad+ \int_0^\tau \int_{\Omega} \Grad \bfphi : {\rm d}\mu_{\mathcal{R}} + \int_0^\tau \intO{ \left< Y_{t,x}; \vr \right> \Grad\Phi \cdot \bfphi } \dt 
			\end{split}
		\end{equation}
		for a.a. $\tau \in (0,T)$ and for any $\bfphi \in C^1([0,T] \times \Ov{\Omega}; R^N)$, $\bfphi\cdot \vc{n}|_{\partial \Omega} = 0$,  where
		\[
		\mu_{\mathcal{R}} \in \mathcal{M}\left([0,T] \times \Ov{\Omega}; R^{N \times N}\right)
		\]
		is a signed (tensor--valued) Radon measure
		on $[0,T] \times \Ov{\Omega}$;
		
		\item
		\begin{equation} \label{MV24}
			\begin{split}
				&\left[
				\intO{ \left< Y_{t,x}; \vr Z \left( s \right) \right> \varphi  } \right]_{t = 0}^{t = \tau} \\ 
				&\geq \int_0^\tau \intO{ \Big[ \left< Y_{t,x} ; \vr Z\left( s \right) \right>  \partial_t \varphi +
					\left< Y_{t,x} ; Z \left( s \right) \vm \right> \cdot \Grad \varphi  \Big] } \dt
			\end{split}
		\end{equation}
		for a.a. $\tau \in (0,T)$, any $\varphi \in C^\infty([0,T] \times R^N)$, $\varphi \geq 0$, $Z' \geq 0$, $Z''\leq 0$, $Z \leq Z_\infty < \infty$;
		
		\item
		\begin{equation} \label{MV25}
			\left[  \intO{ \left< Y_{t,x} ;  \frac{1}{2} \frac{|\vm|^2 }{\vr} + E_{\rm int}  \right> } \right]_{t=0}^{t = \tau} +
			\mathcal{D}(\tau) = \int_0^\tau \intO{ \left< Y_{t,x}; \vm \right> \cdot \Grad\Phi } \dt
		\end{equation}
		where the dissipation defect $\mathcal{D}$ dominates the measure $\mu_{\mathcal{R}}$,
		\begin{equation} \label{MV27}
			\| \mu_{\mathcal{R}}
			\|_{\mathcal{M}([0,\tau) \times \Omega; R^{N \times N})} \leq c \int_0^\tau \mathcal{D}(t) \ \dt
		\end{equation}
		for a.a. $\tau \in (0,T)$. 
	\end{itemize}
\end{Definition}

The existence of measure--valued solutions for given initial data and global in time can be easily shown, for instance by the method of artificial viscosity.
They can be also identified as limits of certain numerical schemes, cf. \cite{FjoMisTad12,FjoMisTad16}.

\subsection{Relative energy}
\label{R}

The \emph{dissipative} measure--valued solutions enjoy certain stability properties, in particular within the class of smooth solutions - the weak--strong uniqueness property shown in \cite{BreFei17}. The key tool is the relative energy inequality.

We start by introducing the relative energy
written in the variables $[\vr, E_{\rm int}, \vm]$ as
\begin{equation} \label{R1}
	\begin{split}
		&\mathcal{E}_Z \left(\vr, E_{\rm int}, \vm \ \Big| r, \Theta, \vc{U} \right) \\
		&=\frac{1}{2} \vr \left| \frac{\vm}{\vr} - \vc{U} \right|^2 + E_{\rm int} - \Theta \vr Z(s(\vr, E_{\rm int})) -
		\frac{ \partial H_\Theta(r, \Theta) }{\partial \vr} (\vr - r) - H_\Theta (r, \Theta),
	\end{split}
\end{equation}
where $H_\Theta$ is the ballistic free energy,
\[
H_\Theta (\vr, \vt) = \vr e(\vr, \vt) - \Theta \vr s(\vr, \vt).
\]
It is worth noting that the relative energy (with $Z(s) = s$) coincides, up to a multiplicative factor $\Theta$, with the relative entropy functional introduced by Dafermos \cite{Dafermos79}, see \cite{Feireisl12}.

The important tool in the analysis of the Euler system is the relative energy inequality proved in \cite{BreFei17}:

\begin{equation} \label{R2}
	\begin{split}
		&\left[ \intO{ \left< Y_{t,x} ; \mathcal{E}_Z \left(\vr, E_{\rm int}, \vm \ \Big| r, \Theta, \vc{U} \right) \right> } \right]_{t = 0}^{t = \tau}
		+ \mathcal{D}(\tau) \\
		&\leq
		- \int_0^\tau \intO{ \left[ \left< Y_{t,x} ; \vr Z\left( s(\vr,E_{\rm int}) \right) \right>  \partial_t \Theta +
			\left< Y_{t,x} ;  Z \left( s(\vr,E_{\rm int})  \right) \vm \right> \cdot \Grad \Theta  \right] } \dt\\
		&\quad+ \int_0^\tau \intO{ \left[ \left< Y_{t,x}; \vr \right> s(r, \Theta) \partial_t \Theta + \left< Y_{t,x}; \vm \right> \cdot s(r, \Theta) \Grad \Theta    \right] } \dt
		\\
		&\quad+ \int_0^\tau \int_\Omega \bigg[ \left< Y_{t,x}; \vr \vc{U} - \vm \right> \cdot \partial_t \vc{U} +
		\left< Y_{t,x}; (\vr \vc{U} - \vm ) \otimes \frac{\vm}{\vr} \right> : \Grad \vc{U} \\
		&\qquad\qquad\qquad- \left< Y_{t,x}; p(\vr, E_{\rm int}) \right> \Div \vc{U}  \bigg] \dx \dt \\
		&\quad + \int_0^\tau \intO{ \left[ \left< Y_{t,x} ; r - \vr \right> \frac{1}{r} \partial_t p(r, \Theta) -
			\left< Y_{t,x} ; \vm \right> \cdot \frac{1}{r} \Grad p(r, \Theta)    \right] } \ \dt\\
		&\quad + \int_0^\tau \intO{ \Grad \Phi \cdot \left< Y_{t,x}; \vm - \vr \vc{U} \right> } \dt - \int_0^\tau \Grad \vc{U}:{\rm d}\mu_{\mathcal{R}}
	\end{split}
\end{equation}
for any trio of test functions $[r, \Theta, \vc{U}]$ belonging to the class
\[
r > 0, \ \Theta > 0, \ r, \Theta \in C^1([0,T]  \times \Ov{\Omega}), \ \vc{U} \in C^1([0,T] \times \Ov{\Omega}; R^N), \
\vc{U} \cdot \vc{n}|_{\partial \Omega} = 0.
\]

\begin{Remark} \label{NR5}
	As a matter of fact, relation (\ref{R2}) was proved in \cite{BreFei17} for the periodic boundary conditions and for the homogeneous case $\Phi=0$. Adaptation to other kinds of boundary conditions like impermeable boundary and to driven fluids, i.e. $\Phi\neq 0$, is straightforward.
\end{Remark}

\section{A singular limit problem}
\label{S}

We illustrate the strength of the theory considering a singular limit problem
for strongly stratified fluids
arising in meteorology or astrophysics, see Majda \cite{Majda03}.

\subsection{Scaled system}

We consider a scaled system driven by an external gradient type force:
\begin{align}
	\partial_t \vr + \Div \vm &= 0, \label{S1}\\
	\partial_t \vm + \Div \left( \frac{\vm \otimes \vm}{\vr} \right)  + \frac{1}{\ep^2} \Grad p(\vr, E_{\rm int}) &= \frac{1}{\ep^2} \vr \Grad \Phi, \label{S2}\\
	\partial_t \left( \frac{1}{2} \frac{|\vm|^2}{\vr} + \frac{1}{\ep^2} E_{\rm int} \right) + \Div \left[ \left( \frac{1}{2} \frac{|\vm|^2}{\vr} + \frac{1}{\ep^2} E_{\rm int}
	+ \frac{1}{\ep^2} p(\vr, E_{\rm int})  \right) \frac{\vm}{\vr} \right] &= \frac{1}{\ep^2} \Grad \Phi \cdot \vm \label{S3}
\end{align}
where $\Phi = \Phi(\vc{x})$ is a given potential.

For the sake of simplicity, we focus on the physically relevant case, where the thermodynamic functions correspond to the perfect gas,
\[
E_{\rm int} = c_v \vr \vt, \ p = \vr \vt, \ s = \log\left( \frac{\vt^{c_v}}{\vr} \right), \ c_v > 0 -
\mbox{specific heat at constant volume.}
\]

We consider the spatial domain to be an infinite slab,
\[
\Omega = \mathcal{T}^2 \times (0,1), \ \mathcal{T}^2 = [0,1]|_{\{ 0,1 \}}
- \mbox{the two dimensional torus,}
\]
meaning the all quantities are space--periodic (with period 1) with respect to the horizontal variable $\vc{x}_h = (x,y)$. The differential operators acting
only on the horizontal variables will be denoted $\nabla_h$, ${\rm div}_h$, etc. We denote $z$ the vertical variable.
We impose the impermeability condition on the lateral boundary,
\[
\vm (\vc{x}_h,z) \cdot \vc{n} = m^3(\vc{x}_h,z) = 0 \ \mbox{for}\ z = 0,1.
\]
Finally, we assume
\begin{equation} \label{NS2}
	\Phi = \Phi (z) = -z,\ \Grad \Phi = [0,0,-1].
\end{equation}

Our goal is to study the singular limit for $\ep \to 0$ in the isothermal regime characterized by constant temperature. More specifically, we consider
solutions close to the static state $[ \vr_s, \Ov{\Theta} ]$,
\begin{equation} \label{NS1}
	\Grad (\vr_s \Ov{\Theta} ) = \vr_s \Grad \Phi, \ \Ov{\Theta} > 0 \ \mbox{a given constant.}
\end{equation}
In view of the ansatz (\ref{NS2}), the static problem (\ref{NS1}) admits a solution determined uniquely by its total mass. Accordingly, there exists a constant $c_0>0$ such that
\[
\vr_s = c_0 \,\exp \left(- \frac{z}{\Ov{\Theta}} \right).
\]

\subsection{Identifying the asymptotic limit, the main result}

Plugging $\vr = \vr_s$, $\vt = \Ov{\Theta}$ in (\ref{S1}) and the entropy inequality (\ref{Ni11}) we obtain
\[
\Div \vm = 0,\ \Div (s(\vr_s, \Ov{\Theta}) \vm ) \geq 0.
\]
In addition, in view of the boundary conditions satisfied by $\vm$, the second inequality reduces to equality. Using the specific form of
$\vr_s$ we easily deduce
\[
\Grad s(\vr_s, \Ov{\Theta}) \cdot \vm = -\frac{\Grad \Phi(z)\cdot \vm}{\Ov{\Theta}} = \frac{m^3}{\Ov{\Theta}} = 0,
\]
meaning the limit momentum (velocity) possesses only the horizontal component, while
\[
{\rm div}_h \vm = 0.
\]

Accordingly, the limit velocity field $\vc{U} = [U^1, U^2,0]=:[\vc{U}_h,0]$ can be taken as the unique solution
of the 2-d \emph{incompressible} Euler system
\begin{equation} \label{NS3}
	\partial_t \vc{U}_h + \vc{U}_h \cdot \nabla_h \vc{U}_h + \nabla_h \Pi = 0, \ {\rm div}_h \vc{U}_h = 0, \ \vc{x}_h \in \mathcal{T}^2,
\end{equation}
supplemented be the initial data
\begin{equation} \label{NS4}
	\vc{U}(0, \vc{x}) = \vc{U}_0 (\vc{x}_h, z) = [U^1_0 (\vc{x}_h,z), U^2_0(\vc{x}_h,z), 0]=[\vc{U}_{h,0} (\vc{x}_h,z), 0].
\end{equation}

The following result is standard, see e.g. Kato \cite{Kato67}.

\begin{Proposition} \label{NP1}
	Let $\vc{U}_{h,0} = \vc{U}_{h,0} (\vc{x}_h)$ be given,
	\[
	\vc{U}_{h,0} \in W^{k,2}(\mathcal{T}^2; R^2), \ {\rm div}_h \vc{U}_{h,0} = 0, \ k > 2.
	\]
	
	Then problem (\ref{NS3}), (\ref{NS4}) admits a (strong) solution $[\vc{U}_h, \Pi]$ unique in the class
	\begin{equation*}
		\begin{split}
			&\vc{U}_h \in C([0,T]; W^{k,2}(\mathcal{T}^2; R^2)), \partial_t \vc{U}_h, \partial_t \Pi, \nabla_h \Pi \in C([0,T]; W^{k-1}(\mathcal{T}^2; R^2)),\
			\int_{\mathcal{T}^2} \Pi \ \dx_h = 0.
		\end{split}
	\end{equation*}
	
\end{Proposition}

We are ready to formulate the main result concerning the singular limit in the system (\ref{S1}--\ref{S3}).

\begin{Theorem} \label{NTS1}
	
	Let the initial data be given such that
	\[
	\vr_{0,\ep} = \vr_s + \ep \vr^{(1)}_{0,\ep},\
	\vt_{0,\ep} = \Ov{\Theta} + \ep \vt^{(1)}_{0,\ep},\
	\vu_{0, \ep} ,
	\]
	where
	\[
	\begin{split}
	\| \vr^{(1)}_{0,\ep} \|_{L^\infty(\Omega)} + \| \vt^{(1)}_{0,\ep} \|_{L^\infty(\Omega)} + \| \vu_{0,\ep}\|_{L^\infty(\Omega; R^N)} \leq c,\\
	\vr^{(1)}_{0,\ep} \to 0,\ \vt^{(1)}_{0,\ep} \to 0,\ \vu_{0,\ep} \to \vc{U}_0 \ \mbox{in}\ L^1(\Omega) \ \mbox{as}\ \ep \to 0,
	\end{split}
	\]
	and where
	\[
	\vc{U}_0 \in W^{k,2}(\Omega; R^3),\ k > 3, \ \vc{U}_0 = [U^1_0, U^2_0,0], \ {\rm div}_h \vc{U}_0 = 0.
	\]
	
	Let $\{ Y^\ep_{t,x} \}_{(t,x) \in (0,T) \times \Omega}$, $\mathcal{D}^\ep$ be a family of dissipative measure--valued solutions to the
	scaled system (\ref{S1}--\ref{S3}), with the initial data
	\[
	Y^\ep_{0,x} = \delta_{\vr_{0,\ep}, \vr_{0,\ep}\vu_{0, \ep}, c_v \vr_{0,\ep} \vt_{0,\ep} },
	\]
	and satisfying the compatibility condition (\ref{MV27}) with a constant $c$ independent of $\ep$.
	
	Then
	\[
	\mathcal{D}^\ep \to 0 \ \mbox{in}\ L^\infty (0,T),
	\]
	and
	\[
	Y^\ep \to \delta_{\vr_s, \vr_s \vc{U}, c_v \vr_s \Ov{\Theta} }
	\ \mbox{in} \ L^\infty(0,T; \mathcal{M}^+ (\mathcal{F})_{{\rm weak-(*)}}),
	\]
	where $[\vr_s, \Ov{\Theta}]$ is the static state and $\vc{U}$ is the unique solution to the Euler system (\ref{NS3}), (\ref{NS4}).

\end{Theorem}

The rest of the section is devoted to the proof of Theorem \ref{NTS1}.
The initial data considered in Theorem \ref{NTS1} are \emph{well--prepared}, meaning adapted to the limit system. In particular, the effect of acoustic
waves is eliminated.
Note that the limit problem is rather different from the isentropic case studied in
\cite{FKNZ16}.

\subsection{Rescaled relative energy}

Definition \ref{DMV1} can be easily adapted to the scaled driven system (\ref{S1}--\ref{S3}).
The relevant relative energy functional reads
\begin{equation} \label{S4}
	\begin{split}
		&\mathcal{E}_{\ep,Z} \left(\vr, E_{\rm int}, \vm \ \Big| r, \Theta, \vc{U} \right)
		\\
		&=\frac{1}{2} \vr \left| \frac{\vm}{\vr} - \vc{U} \right|^2 + \frac{1}{\ep^2} \left[ E_{\rm int} - \Theta \vr Z(s(\vr, E_{\rm int})) -
		\frac{ \partial H_\Theta(r, \Theta) }{\partial \vr} (\vr - r) - H_\Theta (r, \Theta) \right].
	\end{split}
\end{equation}
along with the relative energy inequality
\begin{equation} \label{S5}
	\begin{split}
		&\left[ \intO{ \left< Y_{t,x} ; \mathcal{E}_{\ep,Z} \left(\vr, E_{\rm int}, \vm \ \Big| r, \Theta, \vc{U} \right) \right> } \right]_{t = 0}^{t = \tau}
		+ \mathcal{D}^\ep (\tau)
		\\
		&\leq
		- \frac{1}{\ep^2} \int_0^\tau \intO{ \left[ \left< Y_{t,x} ; \vr Z\left( s(\vr,E_{\rm int}) \right) \right>  \partial_t \Theta +
			\left< Y_{t,x} ;  Z \left( s(\vr,E_{\rm int})  \right) \vm \right> \cdot \Grad \Theta  \right] } \dt\\
		&\quad+ \frac{1}{\ep^2} \int_0^\tau \intO{ \left[ \left< Y_{t,x}; \vr \right> s(r, \Theta) \partial_t \Theta + \left< Y_{t,x}; \vm \right> \cdot s(r, \Theta) \Grad \Theta    \right] } \dt
		\\
		&\quad+ \int_0^\tau \int_\Omega \bigg[ \left< Y_{t,x}; \vr \vc{U} - \vm \right> \cdot \partial_t \vc{U} +
		\left< Y_{t,x}; ( \vr \vc{U} - \vm ) \otimes \frac{\vm}{\vr} \right> : \Grad \vc{U} \\
		&\qquad\qquad\qquad- \frac{1}{\ep^2} \left< Y_{t,x}; p(\vr, E_{\rm int}) \right> \Div \vc{U}  \bigg] \dx \dt \\
		&\quad+ \frac{1}{\ep^2} \int_0^\tau \intO{ \left[ \left< Y_{t,x} ; r - \vr \right> \frac{1}{r} \partial_t p(r, \Theta) -
			\left< Y_{t,x} ; \vm \right> \cdot \frac{1}{r} \Grad p(r, \Theta)    \right] } \dt\\
		&\quad+ \frac{1}{\ep^2} \int_0^\tau \intO{ \Grad \Phi \cdot \left< Y_{t,x}; \vm - \vr \vc{U} \right> } \dt - \int_0^\tau \Grad \vc{U}:{\rm d}\mu_{\mathcal{R}_\ep}.
	\end{split}
\end{equation}

\subsection{Uniform bounds}

Our goal is to establish uniform bounds for the family $Y^\ep$, $\mathcal{D}^\ep$ of measure--valued solutions satisfying the hypotheses of
Theorem \ref{NTS1}.
For $\Theta = \Ov{\Theta} > 0$, $r = \vr_s$ the solution of the stationary problem (\ref{NS1}), the relative energy inequality simplifies to
\[
\begin{split}
&\intO{ \left< Y^\ep_{\tau,x} ; \mathcal{E}_{\ep,Z} \left(\vr, E_{\rm int}, \vm \ \Big| \vr_s, \Ov{\Theta}, \vc{U} \right) \right> }
+ \mathcal{D}^\ep (\tau)
\\
&\leq \int_0^\tau \int_\Omega \bigg[ \left< Y^\ep_{t,x}; \vr \vc{U} - \vm \right> \cdot \partial_t \vc{U} +
\left< Y^\ep_{t,x}; (\vr \vc{U} - \vm ) \otimes \frac{\vm}{\vr} \right> : \Grad \vc{U} \\
&\qquad\qquad\qquad- \frac{1}{\ep^2} \left< Y^\ep_{t,x}; p(\vr, E_{\rm int}) \right> \Div \vc{U}  \bigg] \dx \dt \\
&\quad- \frac{1}{\ep^2} \int_0^\tau \intO{ \left< Y^\ep_{t,x}; \vr \right> \Grad \Phi \cdot \vc{U} } \dt  - \int_0^\tau \Grad \vc{U}:{\rm d}\mu_{\mathcal{R}_\ep} \\
&\quad+ \intO{ \mathcal{E}_{\ep,Z} \left(\vr_{0,\ep}, c_v \vr_{0,\ep} \vt_{0,\ep} , \vr_{0,\ep} \vu_{0,\ep} \ \Big|\  \vr_s , \Ov{\Theta}, \vc{U}(0, \cdot) \right)}.
\end{split}
\]
Thanks to our choice of the initial data, the right--hand side of the above inequality is bounded independently of the cut--off function $Z$, and we get
\begin{equation} \label{S8}
	\begin{split}
		&\intO{ \left< Y^\ep_{\tau,x} ; \mathcal{E}_{\ep} \left(\vr, E_{\rm int}, \vm \ \Big| \vr_s, \Ov{\Theta}, \vc{U} \right) \right> }
		+ \mathcal{D}^\ep (\tau)
		\\
		&\leq \int_0^\tau \int_\Omega \bigg[ \left< Y^\ep_{t,x}; \vr \vc{U} - \vm \right> \cdot \partial_t \vc{U} +
		\left< Y^\ep_{t,x}; (\vr \vc{U} - \vm ) \otimes \frac{\vm}{\vr} \right> : \Grad \vc{U} \\
		&\qquad\qquad\qquad- \frac{1}{\ep^2} \left< Y^\ep_{t,x}; p(\vr, E_{\rm int}) \right> \Div \vc{U}  \bigg] \dx \dt \\
		&\quad- \frac{1}{\ep^2} \int_0^\tau \intO{ \left< Y^\ep_{t,x}; \vr \right> \Grad \Phi \cdot \vc{U} } \dt  - \int_0^\tau \Grad \vc{U}:{\rm d}\mu_{\mathcal{R}_\ep} \\
		&\quad+ \intO{ \mathcal{E}_{\ep} \left(\vr_{0,\ep}, c_v \vr_{0,\ep} \vt_{0,\ep} , \vr_{0,\ep} \vu_{0,\ep} \ \Big|\  \vr_s , \Ov{\Theta}, \vc{U}(0, \cdot) \right)}.
	\end{split}
\end{equation}
where
\[
\begin{split}
&\mathcal{E}_{\ep} \left(\vr, E_{\rm int}, \vm \ \Big| r, \Theta, \vc{U} \right)
\\
&=\frac{1}{2} \vr \left| \frac{\vm}{\vr} - \vc{U} \right|^2 + \frac{1}{\ep^2} \left[ E_{\rm int} - \Theta \vr s(\vr, E_{\rm int}) -
\frac{ \partial H_\Theta(r, \Theta) }{\partial \vr} (\vr - r) - H_\Theta (r, \Theta) \right].
\end{split}
\]

\subsubsection{Entropy estimates}

As the initial data are well--prepared, we have
\[
s(\vr_{0,\ep}, \vt_{0,\ep}) \geq s_0 > - \infty \ \mbox{uniformly for}\ \ep \to 0.
\]
Taking the cut--off function $Z$,
\[
Z'(s) \geq 0,\ Z''(s) \leq 0, Z(s) < 0 \ \mbox{for all} \ s < s_0,\ Z(s) = 0 \ \mbox{whenever}\ s \geq s_0.
\]
we deduce from the entropy inequality (\ref{MV24}) that
\[
\intO{ \left< Y^\ep_{\tau,x}; \vr Z(s) \right> } \geq 0 \ \mbox{for}\ \tau \geq 0.
\]
Consequently,
\[
{\rm supp} [Y^\ep_{\tau,x}] \subset \left\{ [\vr, E_{\rm int}, \vm] \ \Big|\
\vr \geq 0,\ E_{\rm int} \geq 0, \ s(\vr, E_{\rm int}) \geq s_0 \right\},
\]
or, equivalently,
\begin{equation} \label{S9}
	{\rm supp} [Y^\ep_{\tau,x}] \subset \left\{ [\vr, E_{\rm int}, \vm] \ \Big|\
	0 \leq \vr^{1 + \frac{1}{c_v}} \leq c(s_0) E_{\rm int} \right\}\ \mbox{for a.a.}\ (\tau, x).
\end{equation}

\subsubsection{Energy estimates}
\label{EE}

For $\vc{U} = 0$ relation (\ref{S8}) gives rise to
\begin{equation} \label{S10}
	\begin{split}
		&\int_\Omega \bigg< Y^\ep_{\tau,x} ; \frac{1}{2}  \frac{|\vm|^2}{\vr}  + \frac{1}{\ep^2} \bigg[ E_{\rm int} - \Ov{\Theta} \vr s(\vr, E_{\rm int}) \\
		&\qquad\qquad\qquad-\frac{ \partial H_\Theta(\vr_s, \Ov{\Theta}) }{\partial \vr} (\vr_s - r) - H_\Theta (\vr_s, \Ov{\Theta}) \bigg] \bigg> \dx \leq c
	\end{split}
\end{equation}
uniformly for $\ep \to 0$.

Following \cite{BreFei17} we fix a compact set $K \subset (0, \infty)^2$ containing all points
$[\vr_s(x), \Ov{\Theta}]$, $x \in \Omega$ and denote by $\tilde K \subset (0, \infty)^2$ its image in the new phase variables
\[
(\vr, \vt) \mapsto [\vr, c_v \vr \vt] : (0,\infty)^2 \to (0, \infty)^2.
\]
Next, we consider a function $\Psi = \Psi (\vr, E_{\rm int})$,
\[
\Psi \in \DC (0,\infty)^2, \ 0 \leq \Psi \leq 1, \Psi|_{\mathcal{U}} = 1, \ \mbox{where}\ \mathcal{U} \ \mbox{is an open neighborhood of} \ \tilde K
\ \mbox{in}\ (0,\infty)^2.
\]
Finally,
for a measurable function $G(\vr, E_{\rm int}, \vm)$, we set
\[
G = G_{\rm ess} + G_{{\rm res}} , \ G_{\rm ess} = \Psi (\vr, E_{\rm int}) G(\vr, E_{\rm int}, \vm),\
G_{\rm res} = (1 - \Psi (\vr, E_{\rm int})) G(\vr, E_{\rm int}, \vm).
\]
We report the estimate
\cite[Chapter 3, Proposition 3.2]{FeiNov09}
\begin{equation} \label{S11}
	\begin{split}
		&\frac{1}{2}  \frac{|\vm|^2}{\vr}  + \frac{1}{\ep^2} \left[ E_{\rm int} - \Ov{\Theta} \vr s(\vr, E_{\rm int}) -
		\frac{ \partial H_\Theta(\vr_s, \Ov{\Theta}) }{\partial \vr} (\vr_s - r) - H_\Theta (\vr_s, \Ov{\Theta}) \right]
		\\ &\geq c \left(
		\left| \frac{\vm}{\vr}\right|^2 +
		\left[ \left| \frac{\vr - \vr_s}{\ep} \right|^2 + \left| \frac{E_{\rm int} - c_v \vr_s \Ov{\Theta}}{\ep}  \right|^2  \right]_{\rm ess} +
		\left[ \frac{1 + \vr + \vr |s(\vr,E_{\rm int})| + E_{{\rm int}} }{\ep^2}  \right]_{\rm res} \right),
	\end{split}
\end{equation}
where $c$ is a structural constant independent of $\ep$.

Thus relation (\ref{S10}) gives rise to
\begin{equation} \label{S12}
	\intO{ \left< Y^\ep_{\tau,x} ; \frac{|\vm|^2 }{\vr} \right> } \leq c
\end{equation}
and
\begin{align} \label{S13}
	\intO{ \left< Y^\ep_{\tau,x} ; \left[ \vr - \vr_s \right]_{\rm ess}^2 + \left| \left[ \vr - \vr_s \right]_{\rm res} \right| \right> } \leq \ep^2 c , \\
	\label{S14}
	\intO{ \left< Y^\ep_{\tau,x} ; \left[ E_{\rm int} - c_v \vr_s \Ov{\Theta} \right]_{\rm ess}^2 + \left| \left[ E_{\rm int} - c_v \vr_s
		\Ov{\Theta} \right]_{\rm res} \right| \right> } \leq \ep^2 c,
\end{align}
uniformly in $\tau \in [0,T]$, $\ep \to 0$.

\subsubsection{Momentum estimates}
\label{ME}

In view of the estimates established in Section \ref{EE},
we may suppose that
\begin{equation} \label{Mlim}
	\xi^\ep = \left< Y^\ep_{t,x}; \vm \right> \to \xi \ \mbox{weakly-(*) in}\ L^\infty(0,T; L^\alpha(\Omega)),\ \mbox{for some}\ \alpha > 1,
\end{equation}
passing to a subsequence as the case may be.

Indeed we have
\begin{equation*} 
	\left< Y^\ep_{t,x}; \vm \right> = \left< Y^\ep_{t,x}; \sqrt{\vr} \frac{|\vm|}{\sqrt{\vr}} \right>.
\end{equation*}
By Jensen's inequality
\[
\left| \left< Y^\ep_{t,x}; \vm \right> \right|^\alpha \leq \left< Y^\ep_{t,x}; \vr^{\alpha/2} \left( \frac{|\vm|}{\sqrt{\vr}} \right)^\alpha \right>
\lesssim \left< Y^\ep_{t,x}; \frac{|\vm|^2}{{\vr}}  \right> + \left< Y^\ep_{t,x}; \vr^{\alpha/(2 - \alpha)}  \right>.
\]
Thus taking $\alpha > 1$ close to $1$ we may use (\ref{S9}), (\ref{S12}) to obtain the desired conclusion.

\subsection{The limit for $\ep \to 0$}

Let $Y^\ep_{t,x}$, $\mathcal{D}^\ep$ be a family of measure--valued solutions satisfying the hypotheses of Theorem \ref{NTS1}.
In view of the coercivity properties (\ref{S11}), the
conclusion of Theorem \ref{NTS1} follows as soon as we show that
\[
\begin{split}
&\intO{ \left< Y^\ep_{\tau,x} ; \mathcal{E}_{\ep} \left(\vr, E_{\rm int}, \vm \ \Big| \vr_s, \Ov{\Theta}, \vc{U} \right) \right> }
+ \mathcal{D}^\ep (\tau) \to 0 \\
&\qquad\qquad\qquad\mbox{as} \ \ep \to 0
\ \mbox{uniformly for a.a.}\ \tau \in (0,T).
\end{split}
\]
To see this, we apply a Gronwall type argument to inequality (\ref{S8}), where we take $\vc{U}$ the solution of the limit problem (\ref{NS3}), (\ref{NS4}).

\subsubsection{Momentum limit}

With (\ref{Mlim}) in mind, we look at (\ref{MV22}), which is equivalent to
\begin{equation*}
	\left[ \intO{
		\left< Y_{t,x}^\ep; \vr -\vr_s\right> \varphi } \right]_{t = 0}^{t = \tau} =
	\int_0^\tau \intO{ \left[ \left< Y_{t,x}^\ep; \vr - \vr_s \right> \partial_t \varphi + \left< Y_{t,x}^\ep; \vm \right> \cdot \Grad \varphi \right] } \dt
\end{equation*}
because
\begin{equation*}
	\left[ \intO{
		\vr_s \varphi } \right]_{t = 0}^{t = \tau} =
	\int_0^\tau \intO{ \vr_s \partial_t \varphi } \dt.
\end{equation*}
Using the uniform estimates (\ref{S13}) we obtain in the limit $\ep \to 0$:
\begin{equation*}
	\int_0^\tau \intO{ \xi \cdot \Grad \varphi } \dt = 0
\end{equation*}
or, equivalently, 
\begin{equation} \label{S16}
	\Div \xi = 0,\ \xi \cdot \vc{n}|_{\partial \Omega} = 0\quad\mbox{ in the sense of distributions.}
\end{equation}

Similarly, using (\ref{S13}), (\ref{S14}), we may pass to the limit in the entropy balance (\ref{MV24}):
\[
\Div \left( s(\vr_s, \Ov{\Theta}) \xi \right) \geq 0
\]
but since $\xi$ has zero normal trace, this reduces to
\begin{equation} \label{S17}
	\Div \left( s(\vr_s, \Ov{\Theta}) \xi \right) = 0.
\end{equation}
Of course all these relations are understood in the sense of distributions.

Relations (\ref{S16}), (\ref{S17}) are compatible only if
\[
\Div \xi = 0,\ \Grad \vr_s \cdot \xi = 0 \ \mbox{a.e. in}\ (0,T) \times \Omega,
\]
in other words,
\begin{equation} \label{S18}
	\xi^3 = 0,\ {\rm div}_h \xi = 0\ \mbox{a.e. in}\ (0,T) \times \Omega.
\end{equation}

\subsubsection{Geometry and the limit problem}

For $\vc{U}$ - the solution of the limit problem - 
inequality (\ref{S8}) simplifies considerably yielding
\begin{equation} \label{S20}
	\begin{split}
		&\intO{ \left< Y^\ep_{\tau,x} ; \mathcal{E}_{\ep} \left(\vr, E_{\rm int}, \vm \ \Big| \vr_s, \Ov{\Theta}, \vc{U} \right) \right> }
		+ \mathcal{D}^\ep (\tau)
		\\
		&\leq \int_0^\tau \intO{ \left[ \left< Y^\ep_{t,x}; \vr \vc{U} - \vm \right> \cdot \partial_t \vc{U} +
			\left< Y^\ep_{t,x}; ( \vr \vc{U} - \vm ) \otimes \frac{\vm}{\vr} \right> : \Grad \vc{U}  \right] } \dt \\
		&\quad+ c \int_0^\tau \left[ \intO{ \left< Y^\ep_{\tau,x} ; \mathcal{E}_{\ep} \left(\vr, E_{\rm int}, \vm \ \Big| \vr_s, \Ov{\Theta}, \vc{U} \right) \right> }
		+ \mathcal{D}^\ep (\tau) \right] \dt
		+ \omega(\ep),
	\end{split}
\end{equation}
where $\omega = \omega(\ep)$ denotes a generic function enjoying the property $\omega(\ep) \to 0$ as $\ep \to 0$.

Next, we may replace $\frac{\vm}{\vr}$ by $\vc{U}$ keeping (\ref{S20}) still valid. Moreover, in accordance with (\ref{S13}), we may also replace $\vr$ by $\vr_s$ thus obtaining
\begin{equation} \label{S21}
	\begin{split}
		&\intO{ \left< Y^\ep_{\tau,x} ; \mathcal{E}_{\ep} \left(\vr, E_{\rm int}, \vm \ \Big| \vr_s, \Ov{\Theta}, \vc{U} \right) \right> }
		+ \mathcal{D}^\ep (\tau)
		\\
		&\leq \int_0^\tau \intO{ \left( \vr_s \vc{U} - \left< Y^\ep_{t,x}; \vm \right>  \right) \left( \partial_t \vc{U} + \vc{U} \cdot \Grad \vc{U} \right) } \dt \\
		&\quad+ c \int_0^\tau \left[ \intO{ \left< Y^\ep_{\tau,x} ; \mathcal{E}_{\ep} \left(\vr, E_{\rm int}, \vm \ \Big| \vr_s, \Ov{\Theta}, \vc{U} \right) \right> }
		+ \mathcal{D}^\ep (\tau) \right] \dt
		+ \omega(\ep).
	\end{split}
\end{equation}

As $\vc{U}$ satisfies (\ref{NS3}),
relation (\ref{S21}) reduces to
\begin{equation} \label{S24}
	\begin{split}
		&\intO{ \left< Y^\ep_{\tau,x} ; \mathcal{E}_{\ep} \left(\vr, E_{\rm int}, \vm \ \Big| \vr_s, \Ov{\Theta}, \vc{U} \right) \right> }
		+ \mathcal{D}^\ep (\tau)
		\\
		&\lesssim \int_0^\tau \intO{ \left( \xi - \vr_s \vc{U}   \right) \cdot \nabla_h \Pi } \dt \\
		& \quad+ \int_0^\tau \left[ \intO{ \left< Y^\ep_{\tau,x} ; \mathcal{E}_{\ep} \left(\vr, E_{\rm int}, \vm \ \Big| \vr_s, \Ov{\Theta}, \vc{U} \right) \right> }
		+ \mathcal{D}^\ep (\tau) \right] \dt
		+ \omega(\ep).
	\end{split}
\end{equation}

Finally, thanks to (\ref{NS3}), (\ref{S18}),
\[
\intO{ \left( \xi - \vr_s \vc{U}   \right) \cdot \nabla_h \Pi } = 0.
\]
Thus the desired convergence
\[
\intO{ \left< Y^\ep_{\tau,x} ; \mathcal{E}_{\ep} \left(\vr, E_{\rm int}, \vm \ \Big| \vr_s, \Ov{\Theta}, \vc{U} \right) \right> }
+ \mathcal{D}^\ep (\tau) \to 0 \ \mbox{as}\ \ep \to 0 \ \mbox{uniformly in}\ \tau \in (0,T)
\]
follows from Gronwall's lemma.
We have proved Theorem \ref{NTS1}.

\begin{Remark}
	With only minor modifications in the proof, it is also possible to show the same convergence result as in Theorem \ref{NTS1} for the \emph{homogeneous} Euler system in two space dimensions. In this case - since $\Phi=0$ - the limit system is the 2-d incompressible Euler system, too.
\end{Remark}

\def\cprime{$'$} \def\ocirc#1{\ifmmode\setbox0=\hbox{$#1$}\dimen0=\ht0
  \advance\dimen0 by1pt\rlap{\hbox to\wd0{\hss\raise\dimen0
  \hbox{\hskip.2em$\scriptscriptstyle\circ$}\hss}}#1\else {\accent"17 #1}\fi}


\end{document}